\documentclass[12pt,reqno,a4wide]{article}
\usepackage{eurosym}
\usepackage{a4wide}
\usepackage{comment}
\usepackage{mathptmx}
\usepackage{amssymb}
\usepackage{amsfonts}
\usepackage{amsthm}
\usepackage{amsmath}
\usepackage{dsfont}
\usepackage{graphicx}
\usepackage{shadow}
\usepackage[all]{xy}
\usepackage{enumerate}
\usepackage{makeidx}
\usepackage{mathrsfs}
\usepackage{upgreek}
\usepackage{stmaryrd}
\usepackage[utf8]{inputenc}
\usepackage{array}
\usepackage{tikz-cd}
\usepackage{thumbpdf}
\usepackage[colorlinks=true,pagebackref]{hyperref}
\usepackage{textcomp}
\usepackage{chapterbib}
\usepackage{graphics}
\usepackage{longtable}
\usepackage{epsf}
\usepackage{bm}
\usepackage{adjustbox}
\usepackage{centernot}

\makeindex
\newtheorem{theorem}{Theorem}[section]
\newtheorem{lem}[theorem]{Lemma}
\newtheorem{thm}[theorem]{Theorem}
\newtheorem{prop}[theorem]{Proposition}
\newtheorem{cor}[theorem]{Corollary}
\theoremstyle{definition}
\newtheorem*{Beweis}{Proof}
\newtheorem{defn}[theorem]{Definition}
\newtheorem{definition}[theorem]{Definition}

\newtheorem{rem}[theorem]{Remark}

\newtheorem{punto}[theorem]{}

\theoremstyle{remark}

\newtheorem{ex}[theorem]{Example}
\newtheorem{exs}[theorem]{Examples}

\input{tcilatex}
\begin{document}

\title{note on flat semimodules and von Neumann regular semirings\thanks{%
MSC2010: Primary 16Y60; Secondary 16D40, 16E50, 16S50 \newline
Key Words: Semirings; Semimodules; Flat Semimodules; Exact Sequences; Von
Neumann Regular Semirings; Matrix Semirings \newline
The authors would like to acknowledge the support provided by the Deanship
of Scientific Research (DSR) at King Fahd University of Petroleum $\&$
Minerals (KFUPM) for funding this work through projects No. RG1304-1 $\&$
RG1304-2}}
\author{$%
\begin{array}{ccc}
\text{Jawad Abuhlail}\thanks{\text{Corresponding Author}} &  & \text{Rangga
Ganzar Noegraha}\thanks{\text{The paper is extracted from his Ph.D.
dissertation under the supervision of Prof. Jawad Abuhlail.}} \\
\text{abuhlail@kfupm.edu.sa} &  & \text{rangga.gn@universitaspertamina.ac.id}
\\
\text{Department of Mathematics and Statistics} &  & \text{Universitas
Pertamina} \\
\text{King Fahd University of Petroleum $\&$ Minerals} &  & \text{Jl. Teuku
Nyak Arief} \\
\text{31261 Dhahran, KSA} &  & \text{Jakarta 12220, Indonesia}%
\end{array}%
$}
\date{\today }
\maketitle

\begin{abstract}
Flat modules play an important role in the study of the category of modules
over rings and in the characterization of some classes of rings. We study
the $e$\emph{-flatness} for semimodules introduced by the first author using
his new notion of \textit{exact sequences} of semimodules and its
relationships with other notions of flatness for semimodules over semirings.
We also prove that a subtractive semiring over which every right (left)
semimodule is $e$-flat is a von Neumann regular semiring.
\end{abstract}


\section*{Introduction}

\emph{Semirings} are, roughly, rings not necessarily with subtraction. They
generalize both rings and distributive bounded lattices and have, along with
their {\emph{semimodules}, many applications in Computer Science and
Mathematics (e.g., \cite{HW1998, Gla2002, LM2005}). Some applications can be
found in Golan's book \cite{Gol1999}, which is our main reference on this
topic.}

\bigskip

A systematic study of semimodules over semirings was carried out by M.
Takahashi in a series of papers 1981-1990. However, he defined two main
notions in a way that turned out to be \emph{not} natural. Takahashi's \emph{%
tensor products }\cite{Tak1982b} did not satisfy the expected Universal
Property. On the other hand, \emph{Takahashi's exact sequences} of
semimodules \cite{Tak1981} were defined as if this category were exact,
which is \emph{not} the case (in general).

\bigskip

By the beginning of the 21st century, several researchers began to use a
more natural notion of tensor products of semimodules (cf., \cite{Kat2004})
with which the category of semimodules over a commutative semiring is \emph{%
monoidal} rather than \emph{semimonoidal} \cite{Abu2013}. On the other hand,
several notions of exact sequences were introduced (cf., \cite{Pat2003}),
each of which with advantages and disadvantages. One of the most recent
notions is due to Abuhlail \cite{Abu2014} and is based on an intensive study
of the nature of the category of semimodules over a semiring.

\bigskip

In addition to the \emph{categorical notions} of \emph{flat semimodules}
over a semiring, several other notions were considered in the literature,
e.g., the so called $m$-flat\emph{\ semimodules} \cite{Alt2004} (called
\emph{mono-flat} in \cite{Kat2004}). One reason for the interest of such
notions is the phenomenon that, a commutative semiring all of whose
semimodules are flat is a von Neumann regular\emph{\ ring} \cite[Theorem 2.11%
]{Kat2004}. Using a new notion of exact sequences of semimodules over a
semiring, Abuhlail introduced (\cite{Abu2014-SF}) a \emph{homological notion}
of \emph{exactly flat semimodules}, which we call, for short, $e$\emph{-flat}
\emph{semimodules} assuming that an appropriate $\otimes $ functor preserves
short exact sequences.

\bigskip

The paper is divided into three sections.

\bigskip

In Section 1, we collect the basic definitions, examples and preliminaries
used in this paper. Among others, we include the definitions and basic
properties of \emph{exact sequences }as defined by Abuhlail \cite{Abu2014}.

\bigskip

In Section 2, we investigate the $e$-flat semimodules. A \emph{flat}
semimodule is one which is the direct colimit of \emph{finitely presented}
semimodules \cite{Abu2014-SF}. It was proved by Abuhlail \cite[Theorem 3.6]%
{Abu2014-SF} that flat left $S$-semimodules are $e$-flat. We prove in Lemma %
\ref{ret-flat} and Proposition \ref{sum-flat} that the class of $e$-flat
left $S$-semimodules is closed under retracts and direct sums.

\bigskip

In Section 3, we study von Neumann regular semirings. In Theorem \ref%
{sflatvon}, we show that if $S$ is a (left and right) \emph{subtractive}
semiring each of its right semimodules is $S$-$e$-flat, then $S$ is a \emph{%
von Neumann regular semiring}. Conversely, we prove that if $S$ is von
Neumann regular, then every normally $S$-generated right S-semimodule is $S$-%
$m$-flat.

\section{Preliminaries}

\label{prelim}

\qquad In this section, we provide the basic definitions and preliminaries
used in this work. Any notions that are not defined can be found in {our
main reference \cite{Gol1999}. We refer to \cite{Wis1991} for the
foundations of Module and Ring Theory.}

\begin{defn}
(\cite{Gol1999}) A \textbf{semiring} $(S,+,0,\cdot ,1)$ consists of a
commutative monoid $(S,+,0)$ and a monoid $(S,\cdot ,1)$ such that $0\neq 1$
and%
\begin{eqnarray*}
a\cdot 0 &=&0=0\cdot a\text{ for all }a\in S; \\
a(b+c) &=&ab+ac\text{ and }(a+b)c=ac+bc\text{ for all }a,b,c\in S.
\end{eqnarray*}%
If, moreover, the monoid $(S,\cdot ,1)$ is commutative, then we say that $S$
is a \emph{commutative semiring}. We say that $S$ is \emph{additively
idempotent}, if $s+s=s$ for every $s\in S.$
\end{defn}

\begin{exs}
({\cite{Gol1999})}

\begin{itemize}
\item Every ring is a semiring.

\item Any \emph{distributive bounded lattice} $\mathcal{L}=(L,\vee ,0,\wedge
,1)$ is a commutative semiring.

\item Let $R$ be any ring. The set $\mathcal{I}=(Ideal(R),+,0,\cdot ,R)$ of
(two-sided) ideals of $R$ is a semiring with the usual addition and
multiplication of ideals.

\item The set $(\mathbb{Z}^{+},+,0,\cdot ,1)$ (resp. $(\mathbb{Q}%
^{+},+,0,\cdot ,1),$ $(\mathbb{R}^{+},+,0,\cdot ,1)$) of non-negative
integers (resp. non-negative rational numbers, non-negative real numbers) is
a commutative semiring (resp. \emph{semifield}) which is not a ring (not a
field).

\item $(M_{n}(S),+,0,\cdot ,I_{n}),$ the set of all $n\times n$ matrices
over a semiring $S,$ is a semiring with the usual addition and
multiplication of matrices.

\item $\mathbb{B}:=\{0,1\}$ with $1+1=1$ is a semiring, called the Boolean
semiring.

\item The \emph{max-min algebra} $\mathbb{R}_{\max ,\min }:=(\mathbb{R}\cup
\{-\infty ,\infty \},\max ,-\infty ,\min ,\infty )$ is an additively
idempotent semiring$.$

\item The \emph{log algebra} $(\mathbb{R}\cup \{-\infty ,\infty \},\oplus
,\infty ,+,0)$ is a semiring, where%
\begin{equation*}
x\oplus y=-ln(e^{-x}+e^{-y})
\end{equation*}
\end{itemize}
\end{exs}

\begin{punto}
\cite{Gol1999} Let $S$ and $T$ be semirings. The categories $_{S}\mathbf{SM}$
of \textbf{left} $S$-\textbf{semimodules} with arrows the $S$-linear maps, $%
\mathbf{SM}_{T}$ of right $S$-semimodules with arrows the $T$-linear maps,
and $_{S}\mathbf{SM}_{T}$ of $(S,T)$-bisemimodules are defined in the usual
way (as for modules and bimodules over rings). We write $L\leq _{S}M$ to
indicate that $L$ is an $S$-subsemimodule of the left (right) $S$-semimodule
$M.$
\end{punto}

\begin{ex}
The category of $\mathbb{Z}^{+}$-semimodules is nothing but the category of
commutative monoids.
\end{ex}

\begin{definition}
\cite[page 162]{Gol1999} Let $S$ be a semiring. An equivalence relation $%
\rho $ on a left $S$-semimodule $M$ is a \textbf{congruence relation}, if it
preserves the addition and the scalar multiplication on $M,$ \emph{i.e. }for
all $s\in S$ and $m,m^{\prime },n,n^{\prime }\in M:$%
\begin{equation*}
m\rho m^{\prime }\text{ and }n\rho n^{\prime }\Longrightarrow (m+m^{\prime
})\rho (n+n^{\prime }),
\end{equation*}%
\begin{equation*}
m\rho m^{\prime }\Longrightarrow (sm)\rho (sm^{\prime }).
\end{equation*}
\end{definition}

\begin{punto}
(\cite[page 150, 154]{Gol1999}) Let $S$ be a semiring and $M$ a left $S$%
-semimodule.

\begin{enumerate}
\item The \textbf{subtractive closure }of $L\leq _{S}M$ is defined as%
\begin{equation}
\overline{L}:=\{m\in M\mid \text{ }m+l=l^{\prime }\text{ for some }%
l,l^{\prime }\in L\}.  \label{L-s-closure}
\end{equation}%
We say that $L$ is subtractive if $L=\overline{L}.$ The left $S$-semimodule $%
M$ is a \textbf{subtractive semimodule, }if every $S$-subsemimodule $L\leq
_{S}M$ is subtractive.

\item The set of \textbf{cancellative elements}\emph{\ }of $M$ is defined as%
\begin{equation*}
K^{+}(M)=\{x\in M\mid x+y=x+z\Longrightarrow y=z\text{ for any }y,z\in M\}.
\end{equation*}%
We say that $M$ is a \textbf{cancellative semimodule}, if $K^{+}(M)=M.$
\end{enumerate}
\end{punto}

\begin{punto}
\label{variety}(cf., \cite{AHS2004})\ The category $_{S}\mathbf{SM}$ of left
semimodules over a semiring $S$ is a \emph{variety} (i.e. closed under
homomorphic images, subobjects and arbitrary products), whence complete
(i.e. has all limits, e.g., direct products, equalizers, kernels, pullbacks,
inverse limits) and cocomplete (i.e. has all colimits, e.g., direct
coproducts, coequalizers, cokernels, pushouts, direct colimits).
\end{punto}

\begin{punto}
\label{tensprod}With the tensor product of a right $S$-semimodule $L$ and a
left $S$-semimodule $M$, we mean the commutative monoid $L\otimes _{S}M$ in
the sense of \cite[3.1]{Kat1997}, and not that in the sense of Takahashi
adapted by Golan \cite{Gol1999}, which we denote by $M\boxtimes _{S}N$.
Abuhlail \cite{Abu2013} showed that $M\boxtimes _{S}N=c(M\otimes _{S}N),$
the \emph{cancellative hull }of $M\otimes _{S}N$. See also \cite[2.1]%
{Abu2014-SF}.
\end{punto}

The following results is folklore and is implicit in \cite{Kat1997}.

\begin{lem}
\label{lem258}For every right $S$-semimodule $M$, there exists a natural
right $S$-isomorphism
\begin{equation*}
\theta _{M}:M\otimes _{S}S\rightarrow M,\text{ }m\otimes _{S}s\mapsto ms.
\end{equation*}
\end{lem}

\subsection*{Exact Sequences}

\bigskip

Throughout, $(S,+,0,\cdot ,1)$ is a semiring and, unless otherwise
explicitly mentioned, an $S$-semimodule is a \emph{left }$S$-semimodule.

\bigskip

\begin{definition}
A morphism of left $S$-semimodules $f:L\rightarrow M$ is

$k$-\textbf{normal}, if whenever $f(m)=f(m^{\prime })$ for some $m,m^{\prime
}\in M,$ we have $m+k=m^{\prime }+k^{\prime }$ for some $k,k^{\prime }\in
Ker(f);$

$i$-\textbf{normal}, if $\func{Im}(f)=\overline{f(L)}$ ($:=\{m\in M|\text{ }%
m+l\in L\text{ for some }l\in L\}$).

\textbf{normal}, if $f$ is both $k$-normal and $i$-normal.
\end{definition}

\begin{rem}
Among others, Takahashi (\cite{Tak1981}) and Golan \cite{Gol1999} called $k$%
-normal (resp., $i$-normal, normal) $S$-linear maps $k$\emph{-regular}
(resp., $i$\emph{-regular}, \emph{regular}) morphisms. Our terminology is
consistent with Category Theory noting that the \textbf{normal epimorphisms}
are exactly the normal surjective $S$-linear maps, and the \textbf{normal
monomorphisms} are exactly the normal injective $S$-linear maps (see \cite%
{Abu2014}).
\end{rem}

We often refer to the following technical lemma:

\begin{lem}
\emph{(\cite[Lemma 1.17]{AN2021}) }\label{i-normal}Let $L\overset{f}{%
\rightarrow }M\overset{g}{\rightarrow }N$ be a sequence of semimodules.

\begin{enumerate}
\item Let $g$ be injective.

\begin{enumerate}
\item $f$ is $k$-normal if and only if $g\circ f$ is $k$-normal.

\item If $g\circ f$ is $i$-normal (normal), then $f$ is $i$-normal (normal).

\item Assume that $g$ is $i$-normal. Then $f$ is $i$-normal (normal) if and
only if $g\circ f$ is $i$-normal (normal).
\end{enumerate}

\item Let $f$ be surjective.

\begin{enumerate}
\item $g$ is $i$-normal if and only if $g\circ f$ is $i$-normal.

\item If $g\circ f$ is $k$-normal (normal), then $g$ is $k$-normal (normal).

\item Assume that $f$ is $k$-normal. Then $g$ is $k$-normal (normal) if and
only if $g\circ f$ is $k$-normal (normal).
\end{enumerate}
\end{enumerate}
\end{lem}

The proof of the following lemma is straightforward:

\begin{lem}
\label{u-sum}

\begin{enumerate}
\item Let $\{f_{\lambda }:L_{\lambda }\longrightarrow M_{\lambda
}\}_{\Lambda }$ be a non-empty collection of left $S$-semimodule morphisms
and consider the induced $S$-linear map $f:\bigoplus\limits_{\lambda \in
\Lambda }L_{\lambda }\longrightarrow \bigoplus\limits_{\lambda \in \Lambda
}M_{\lambda }.$ Then $f$ is normal \emph{(}resp. $k$-normal, $i$-normal\emph{%
)} if and only if $f_{\lambda }$ is normal \emph{(}resp. $k$-normal, $i$%
-normal\emph{)} for every $\lambda \in \Lambda .$ 

\item A morphism $\varphi :L\longrightarrow M$ of left $S$-semimodules is
normal \emph{(}resp. $k$-normal, $i$-normal\emph{)} if and only if $\mathrm{%
id}_{F}\otimes _{S}\varphi :F\otimes _{S}L\longrightarrow F\otimes _{S}M$ is
normal \emph{(}resp. $k$-normal, $i$-normal\emph{)} for every non-zero free
right $S$-semimodule $F.$

\item If $P_{S}$ is projective and $\varphi :L\longrightarrow M$ is a normal
\emph{(}resp. $k$-normal, $i$-normal\emph{)} morphism of left $S$%
-semimodules, then $\mathrm{id}_{F}\otimes _{S}\varphi :P\otimes
_{S}L\longrightarrow P\otimes _{S}M$ is normal \emph{(}resp. $k$-normal, $i$%
-normal\emph{)}.
\end{enumerate}
\end{lem}

There are several notions of exactness for sequences of semimodules. In this
paper, we use the relatively new notion introduced by Abuhlail:

\begin{definition}
\label{Abu-exs}(\cite[2.4]{Abu2014}) A sequence
\begin{equation}
L\overset{f}{\longrightarrow }M\overset{g}{\longrightarrow }N  \label{LMN}
\end{equation}%
of left $S$-semimodules is \textbf{exact}, if $g$ is $k$-normal and $%
f(L)=Ker(g).$
\end{definition}

\begin{punto}
\label{def-exact}We call a sequence of $S$-semimodules $L\overset{f}{%
\rightarrow }M\overset{g}{\rightarrow }N$

\emph{proper-exact} if $f(L)=\mathrm{Ker}(g)$ (exact in the sense of
Patchkoria \cite{Pat2003});

\emph{semi-exact} if $\overline{f(L)}=\mathrm{Ker}(g)$ (exact in the sense
of Takahashi \cite{Tak1982a});

\emph{quasi-exact} if $\overline{f(L)}=\mathrm{Ker}(g)$ and $g$ is $k$%
-normal (exact in the sense of Patil and Doere \cite{PD2006}).
\end{punto}

\begin{punto}
We call a (possibly infinite) sequence of $S$-semimodules%
\begin{equation}
\cdots \rightarrow M_{i-1}\overset{f_{i-1}}{\rightarrow }M_{i}\overset{f_{i}}%
{\rightarrow }M_{i+1}\overset{f_{i+1}}{\rightarrow }M_{i+2}\rightarrow \cdots
\label{chain}
\end{equation}%
\emph{chain complex} if $f_{j+1}\circ f_{j}=0$ for every $j;$

\emph{exact} (resp., \emph{proper-exact}, \emph{semi-exact, quasi-exact}) if
each partial sequence with three terms $M_{j}\overset{f_{j}}{\rightarrow }%
M_{j+1}\overset{f_{j+1}}{\rightarrow }M_{j+2}$ is exact (resp.,
proper-exact, semi-exact, quasi-exact).

A \textbf{short exact sequence} (or a \textbf{Takahashi extension} \cite%
{Tak1982b}) of $S$-semimodules is an exact sequence of the form%
\begin{equation*}
0\longrightarrow L\overset{f}{\longrightarrow }M\overset{g}{\longrightarrow }%
N\longrightarrow 0
\end{equation*}
\end{punto}

The following examples show some of the advantages of the new definition of
exact sequences over the old ones:

\begin{lem}
\label{exact}Let $L,M$ and $N$ be $S$-semimodules.

\begin{enumerate}
\item $0\longrightarrow L\overset{f}{\longrightarrow }M$ is exact if and
only if $f$ is injective.

\item $M\overset{g}{\longrightarrow }N\longrightarrow 0$ is exact if and
only if $g$ is surjective.

\item $0\longrightarrow L\overset{f}{\longrightarrow }M\overset{g}{%
\longrightarrow }N$ is proper-exact and $f$ is normal (semi-exact and $f$ is
normal) if and only if $L\simeq \mathrm{Ker}(g).$

\item $0\longrightarrow L\overset{f}{\longrightarrow }M\overset{g}{%
\longrightarrow }N$ is exact if and only if $L\simeq \mathrm{Ker}(g)$ and $g$
is $k$-normal.

\item $L\overset{f}{\longrightarrow }M\overset{g}{\longrightarrow }%
N\longrightarrow 0$ is semi-exact and $g$ is normal if and only if $N\simeq
M/f(L).$

\item $L\overset{f}{\longrightarrow }M\overset{g}{\longrightarrow }%
N\longrightarrow 0$ is exact if and only if $N\simeq M/f(L)$ and $f$ is $i$%
-normal.

\item $0\longrightarrow L\overset{f}{\longrightarrow }M\overset{g}{%
\longrightarrow }N\longrightarrow 0$ is exact if and only if $L\simeq
\mathrm{Ker}(g)$ and $N\simeq M/L.$
\end{enumerate}
\end{lem}

Our definition of exact sequences allows us to recover the following
well-known result for short exact sequences of modules over rings:

\begin{cor}
\label{M/L}The following assertions are equivalent:

\begin{enumerate}
\item $0\rightarrow L\overset{f}{\rightarrow }M\overset{g}{\rightarrow }%
N\rightarrow 0$ is an exact sequence of $S$-semimodules;

\item $L\simeq \mathrm{Ker}(g)$ and $N\simeq Coker(f)$ ($=M/f(L)$);

\item $f$ is injective, $f(L)=\mathrm{Ker}(g),$ $g$ is surjective and ($k$%
-)normal.

In this case, $f$ and $g$ are normal morphisms.
\end{enumerate}
\end{cor}

\begin{rem}
An $S$-linear map is a monomorphism if and only if it is injective. Every
surjective $S$-linear map is an epimorphism. The converse is not true in
general (cf. \cite{WT1989}).
\end{rem}

\begin{lem}
\label{semi-ex}Let%
\begin{equation*}
\xymatrix{A' \ar[r]^{i} \ar[d]_{f} & A \ar[r]^{p} \ar[d]_{g} & A''
\ar@{-->}[d]^{h} \\ B' \ar[r]^{j} & B \ar[r]^{q} & B''}
\end{equation*}%
be a commutative diagram of left $S$-semimodules and semi-exact rows. If $p$
is a normal epimorphism, then there exists a unique $S$-linear map $%
h:A^{\prime \prime }\rightarrow B^{\prime \prime }$ making the augmented
diagram commute.

\begin{enumerate}
\item If, moreover, $q$ is a normal epimorphism, $f$ is surjective and $g$
is injective (an isomorphism), then $h$ is injective (an isomorphism).

\item If, moreover, $A$ and $B$ are cancellative, $j,$ $f$ and $h$ are
injective, then $g$ is injective.
\end{enumerate}
\end{lem}

\begin{Beweis}
Since $p$ is normal, $A^{\prime \prime }\simeq Coker(i)$ by Lemma (\ref%
{exact}) (5). Since $q\circ g\circ i=q\circ j\circ f=0,$ the existence and
uniqueness of $h$ follows directly from the \emph{Universal Property of
Cokernels}. However, we give an elementary proof that $h$ is well-defined
using diagram chasing. Let $a^{\prime \prime }\in A^{\prime \prime }.$ Since
$p$ is surjective, there exists $a\in A$ such that $p(a)=a^{\prime \prime }.$
Consider%
\begin{equation*}
h:A^{\prime \prime }\rightarrow B^{\prime \prime },\text{ }a^{\prime \prime
}\mapsto q(g(a)).
\end{equation*}%
\textbf{Claim:} $h$ is well defined.

Suppose there exist $p(a_{1})=a^{\prime \prime }=p(a_{2}).$ Since $p$ is $k$%
-normal, $a_{1}+k_{1}=a_{2}+k_{2}$ for some $k_{1},k_{2}\in Ker(p)=\func{Im}%
(i)$. Let $a_{1}^{\prime },\widetilde{a}_{1},a_{2}^{\prime },\widetilde{a}%
_{2}\in A^{\prime }$ be such that $k_{1}+i(a_{1}^{\prime })=i(\widetilde{a}%
_{1})$ and $k_{2}+i(a_{2}^{\prime })=i(\widetilde{a}_{2}).$ It follows that%
\begin{equation*}
a_{1}+k_{1}+i(a_{1}^{\prime })=a_{1}+i(\widetilde{a}_{1})\text{ and }%
a_{2}+k_{2}+i(a_{2}^{\prime })=a_{2}+i(\widetilde{a}_{2})
\end{equation*}%
and so%
\begin{equation*}
(q\circ g\circ i)(\widetilde{a}_{1})=(q\circ j\circ f)(\widetilde{a}%
_{1})=0=(q\circ j\circ f)(\widetilde{a}_{2})=(q\circ g\circ i)(\widetilde{a}%
_{2}).
\end{equation*}%
So%
\begin{eqnarray*}
(q\circ g)(a_{1}) &=&(q\circ g)(a_{1})+(q\circ j\circ f)(\widetilde{a}_{1})
\\
&=&(q\circ g)(a_{1})+(q\circ g\circ i)(\widetilde{a}_{1}) \\
&=&(q\circ g)(a+i(\widetilde{a}_{1})) \\
&=&(q\circ g)(a_{1}+k_{1}+i(a_{1}^{\prime })) \\
&=&(q\circ g)(a_{1}+k_{1})+(q\circ j\circ f)(a_{1}^{\prime }) \\
&=&(q\circ g)(a_{2}+k_{2}) \\
&=&(q\circ g)(a_{2}+k_{2})+(q\circ j\circ f)(a_{2}^{\prime }) \\
&=&(q\circ g)(a_{2}+k_{2}+i(a_{2}^{\prime })) \\
&=&(q\circ g)(a_{2}+i(\widetilde{a}_{2})) \\
&=&(q\circ g)(a_{2})+(q\circ j\circ f)(\widetilde{a}_{2}) \\
&=&(q\circ g)(a_{2}).
\end{eqnarray*}%
Thus $h$ is well defined and $h\circ p=q\circ g$ by the definition of $h$.
Clearly, $h$ is unique.

\begin{enumerate}
\item Suppose that $h(x_{1})=h(x_{2})$ for some $x_{1},x_{2}\in A^{\prime
\prime }.$ Since $p$ is surjective, $x_{1}=p(a_{1})$ and $x_{2}=p(a_{2})\ $%
for some $a_{1},$ $a_{2}\in A.$ So,
\begin{equation*}
q(g(a_{1}))=h(p(a_{1}))=h(x_{1})=h(x_{2})=h(p(a_{2}))=q(g(a_{2})).
\end{equation*}%
Since the second row is semi-exact, there exist $y_{1},y_{2}\in \func{Im}(j)$
such that $g(a_{1})+y_{1}=g(a_{2})+y_{2}.$ Let $z_{1},\widetilde{z}%
_{1},z_{2},\widetilde{z}_{2}\in B^{\prime }$ be such that $y_{1}+j(z_{1})=j(%
\widetilde{z}_{1})$ and $y_{2}+j(z_{2})=j(\widetilde{z}_{2}).$ It follows
that%
\begin{eqnarray*}
g(a_{1})+y_{1}+j(z_{1})+j(z_{2}) &=&g(a_{2})+y_{2}+j(z_{2})+j(z_{1}) \\
g(a_{1})+j(\widetilde{z}_{1})+j(z_{2}) &=&g(a_{2})+j(\widetilde{z}%
_{2})+j(z_{1}) \\
g(a_{1})+j(\widetilde{z}_{1}+z_{2}) &=&g(a_{2})+j(\widetilde{z}_{2}+z_{1})
\end{eqnarray*}%
Since $f$ is surjective, there exist $w_{1},w_{2}\in A^{\prime }$ Such that $%
f(w_{1})=\widetilde{z}_{1}+z_{2}$ and $f(w_{2})=\widetilde{z}_{2}+z_{1}.$
So, we have%
\begin{equation*}
\begin{array}{ccccc}
g(a_{1}+i(w_{1})) & = & g(a_{1})+(g\circ i)(w_{1}) & = & g(a_{1})+(j\circ
f)(w_{1}) \\
& = & g(a_{1})+j(\widetilde{z}_{1}+z_{2}) & = & g(a_{2})+j(\widetilde{z}%
_{2}+z_{1}) \\
& = & g(a_{2})+(g\circ i)(w_{2}) & = & g(a_{2}+i(w_{2}))%
\end{array}%
\end{equation*}%
Since $g$ is injective, we have $a_{1}+i(w_{1})=a_{2}+i(w_{2}),$ whence%
\begin{equation*}
x_{1}=p(a_{1})=p(a_{1}+i(w_{1}))=p(a_{2}+i(w_{2}))=p(a_{2})=x_{2}.
\end{equation*}%
It follows that $h$ is injective. If $g$ is surjective, then $h\circ
p=q\circ g$ is surjective, whence $h$ is surjective.

\item Suppose that $g(a_{1})=g(a_{2})$ for some $a_{1},a_{2}\in A.$ It
follows that%
\begin{equation*}
(h\circ p)(a_{1})=(q\circ g)(a_{1})=(q\circ g)(a_{2})=(h\circ p)(a_{2}),
\end{equation*}%
whence $p(a_{1})=p(a_{2})$ ($h$ is injective, by assumption). Since the
first row is semi-exact, there exist $y_{1},y_{2}\in \func{Im}(i)=Ker(p)$
such that $a_{1}+y_{1}=a_{2}+y_{2}.$ Let $w_{1},\widetilde{w}_{1},w_{2},%
\widetilde{w}_{2}\in A^{\prime }$ be such that $y_{1}+i(w_{1})=i(\widetilde{w%
}_{1})$ and $y_{2}+i(w_{2})=i(\widetilde{w}_{2}).$ It follows that%
\begin{equation*}
a_{1}+y_{1}+i(w_{1})=a_{1}+i(\widetilde{w}_{1})\text{ and }%
a_{2}+y_{2}+i(w_{2})=a_{2}+i(\widetilde{w}_{2}).
\end{equation*}%
Consequently, we have%
\begin{eqnarray*}
a_{1}+y_{1}+i(w_{1})+i(w_{2}) &=&a_{2}+y_{2}+i(w_{2})+i(w_{1}) \\
a_{1}+i(\widetilde{w}_{1})+i(w_{2}) &=&a_{2}+i(\widetilde{w}_{2})+i(w_{1}) \\
a_{1}+i(\widetilde{w}_{1}+w_{2})) &=&a_{2}+i(\widetilde{w}_{2}+w_{1})
\end{eqnarray*}%
It follows that%
\begin{eqnarray*}
g(a_{1})+(g\circ i)(\widetilde{w}_{1}+w_{2}) &=&g(a_{2})+(g\circ i)(%
\widetilde{w}_{2}+w_{1}) \\
g(a_{1})+(j\circ f)(\widetilde{w}_{1}+w_{2}) &=&g(a_{2})+(j\circ f)(%
\widetilde{w}_{2}+w_{1}).
\end{eqnarray*}%
Since $B$ is cancellative and both $f$ and $j$ are injective, we conclude
that $\widetilde{w}_{1}+w_{2}=\widetilde{w}_{2}+w_{1}.$ Since $A$ is
cancellative, we conclude $a_{1}=a_{2}.\blacksquare $
\end{enumerate}
\end{Beweis}

\begin{lem}
\label{pull}Consider an exact sequence%
\begin{equation*}
0\longrightarrow L\overset{f}{\longrightarrow }M\overset{g}{\longrightarrow }%
N\longrightarrow 0
\end{equation*}%
of left $S$-semimodules. If $U\leq _{S}N$ is a (subtractive)\ subsemimodule,
then in the pullback $(P,\iota ^{\prime },g^{\prime })$ of $\iota
:U\hookrightarrow N$ and $g:M\longrightarrow N$ the $S$-linear map $%
g^{\prime }:P\rightarrow M$ is a (normal)\ monomorphism.
\end{lem}

\begin{Beweis}
Let $U\leq _{S}N.$ By \cite[1.7]{Tak1982b}, the \emph{pullback} of $\iota
:U\hookrightarrow N$ and $g:M\longrightarrow N$ is $(P,\iota ^{\prime
},g^{\prime })$, where%
\begin{eqnarray}
P &:&=\{(u,m)\in U\times M\text{ }|\text{ }\iota (u)=g(m)\}  \label{pullback}
\\
\iota ^{\prime } &:&P\rightarrow U,\text{ }(u,m)\mapsto u;  \notag \\
g^{\prime } &:&P\rightarrow M,\text{ }(u,m)\mapsto m.  \notag
\end{eqnarray}%
Consider the following diagram of left $S$-semimodules%
\begin{equation}
\xymatrix{& & 0 \ar[d] & 0 \ar[d] & \\ 0 \ar[r] & L \ar[r]^{f'}
\ar@{=}[d]_{id} & P \ar[d]^{g'} \ar[r]^{\iota '} & U \ar[d]^{\iota} \ar[r] &
0 \\ 0 \ar[r] & L \ar[r]_{f} & M \ar[r]_{g} & N \ar[r] & 0}  \label{pbc}
\end{equation}

\textbf{Claim I:}\ $g^{\prime }$ is injective.

Suppose that $g^{\prime }(u_{1},m_{1})=g^{\prime }(u_{2},m_{2})$ for some $%
u_{1},u_{2}\in U$ and $m_{1},m_{2}\in M.$ Then $m_{1}=m_{2}$ and moreover,
\begin{equation*}
\iota (u_{1})=g(m_{1})=g(m_{2})=\iota (u_{2}),
\end{equation*}%
whence $u_{1}=u_{2}$ since $\iota $ is injective. Consequently, $g^{\prime }$
is injective.

\textbf{Claim II:\ }If $\iota (U)\leq _{S}N$ is subtractive, then $g^{\prime
}(P)\leq _{S}M$ is subtractive.

Let $m\in \func{Im}(g^{\prime }),$ i.e. $m+g^{\prime
}(u_{1},m_{1})=g^{\prime }(u_{2},m_{2})$ for some $u_{1},u_{2}\in U$ and $%
m_{1},m_{2}\in M.$ It follows that $m+m_{1}=m_{2},$ whence $%
g(m)+g(m_{1})=g(m_{2})$ and $g(m)+\iota (u_{1})=\iota (u_{2}).$ Since $\iota
$ is a normal monomorphism, it follows that $g(m)=\iota (u)$ for some $u\in
U $ and so $m=g^{\prime }(m,u).$ Consequently, $g^{\prime }(P)=\func{Im}%
(g^{\prime }),$ i.e. $g^{\prime }$ is a normal monomorphism.$\blacksquare $
\end{Beweis}

\begin{prop}
\label{adj-lim}\emph{(}cf., \emph{\cite[Proposition 3.2.2]{Bor1994})} Let $%
\mathfrak{C},\mathfrak{D}$ be arbitrary categories and $\mathfrak{C}\overset{%
F}{\longrightarrow }\mathfrak{D}\overset{G}{\longrightarrow }\mathfrak{C}$
be covariant functors such that $(F,G)$ is an adjoint pair.

\begin{enumerate}
\item $F$ preserves all colimits which turn out to exist in $\mathfrak{C}.$

\item $G$ preserves all limits which turn out to exist in $\mathfrak{D}.$
\end{enumerate}
\end{prop}

\begin{cor}
\label{ad-l-cor}Let $S,$ $T$ be semirings and $_{T}F_{S}$ a $(T,S)$%
-bisemimodule.

\begin{enumerate}
\item $F\otimes _{S}-:$ $_{S}\mathbf{SM}\longrightarrow $ $_{T}\mathbf{SM}$
preserves all colimits.

\item For every non-empty collection of left $S$-semimodules $\{X_{\lambda
}\}_{\Lambda },$ we have a canonical isomorphism of left $T$-semimodules%
\begin{equation*}
F\otimes _{S}\bigoplus\limits_{\lambda \in \Lambda }X_{\lambda }\simeq
\bigoplus\limits_{\lambda \in \Lambda }(F\otimes _{S}X_{\lambda }).
\end{equation*}

\item For any directed system of left $S$-semimodules $(X_{j},\{f_{jj^{%
\prime }}\})_{J},$ we have an isomorphism of left $T$-semimodules%
\begin{equation*}
F\otimes _{S}\lim_{\longrightarrow }X_{j}\simeq \text{ }\lim_{%
\longrightarrow }(F\otimes _{S}X_{j}).
\end{equation*}

\item $F\otimes _{S}-$ preserves coequalizers.

\item $F\otimes _{S}-$ preserves cokernels.
\end{enumerate}
\end{cor}

\begin{Beweis}
The proof can be obtained as a direct consequence of Proposition \ref%
{adj-lim} and the fact that $(F\otimes _{S}-,\mathrm{Hom}_{T}(F,-))$ is an
adjoint pair of covariant functors \cite{KN2011}.$\blacksquare $
\end{Beweis}

\begin{prop}
\label{r-exact}Let $_{T}G_{S}$ be a $(T,S)$-bisemimodule and consider the
functor $G\otimes _{S}-:$ $_{S}\mathbf{SM}\longrightarrow $ $_{T}\mathbf{SM}%
. $ Let%
\begin{equation}
L\overset{f}{\rightarrow }M\overset{g}{\rightarrow }N\rightarrow 0
\label{L-r}
\end{equation}%
be a sequence of left $S$-semimodules and consider the sequence of left $T$%
-semimodules%
\begin{equation}
G\otimes _{S}L\overset{G\otimes f}{\longrightarrow }G\otimes _{S}M\overset{%
G\otimes g}{\longrightarrow }G\otimes _{S}N\rightarrow 0  \label{FL-r}
\end{equation}

\begin{enumerate}
\item If $M\overset{g}{\rightarrow }N\rightarrow 0$ is exact and $g$ is
normal, then $G\otimes _{S}M\overset{G\otimes g}{\longrightarrow }G\otimes
_{S}N\rightarrow 0$ is exact and $G\otimes g$ is normal.

\item If \emph{(\ref{L-r}) }is semi-exact and $g$ is normal, then \emph{(\ref%
{FL-r}) }is semi-exact and $G\otimes g$ is normal.

\item If \emph{(\ref{L-r}) }is exact and $G\otimes _{S}f$ is $i$-normal,
then \emph{(\ref{FL-r}) }is exact.
\end{enumerate}
\end{prop}

\begin{Beweis}
The following implications are obvious: $M\overset{g}{\rightarrow }%
N\rightarrow 0$ is exact $\Longrightarrow $ $g$ is surjective $%
\Longrightarrow $ $G\otimes g$ is surjective $\Longrightarrow $ $G\otimes
_{S}M\overset{G\otimes g}{\longrightarrow }G\otimes _{S}N\rightarrow 0$ is
exact.

\begin{enumerate}
\item Assume that $g$ is normal and consider the exact sequence of $S$%
-semimodules%
\begin{equation*}
0\longrightarrow \mathrm{Ker}(g)\overset{\iota }{\longrightarrow }M\overset{g%
}{\longrightarrow }N\longrightarrow 0.
\end{equation*}%
Then $N\simeq \mathrm{Coker}(\iota ).$ By Corollary \ref{ad-l-cor} (1), $%
G\otimes _{S}-$ preserves cokernels and so $G\otimes g=\mathrm{coker}%
(G\otimes \iota )$ whence normal.

\item Apply Lemma \ref{exact}: The assumptions on (\ref{L-r}) are equivalent
to $N=\mathrm{Co}$\textrm{$k$}$\mathrm{er}(f)$. Since $G\otimes _{S}-$
preserves cokernels, we conclude that $G\otimes _{S}N=\mathrm{Coker}%
(G\otimes f),$ \emph{i.e.} (\ref{FL-r}) is semi-exact and $G\otimes g$ is
normal.

\item This follows directly from (2) and the assumption on $G\otimes
f.\blacksquare $
\end{enumerate}
\end{Beweis}

\section{Flat Semimodules}

\markright{\scriptsize\tt Chapter \ref{sec-flat}: Semirings
and Semimodules} 

The notion of \emph{exactly flat semimodules} was introduced by Abuhlail
\cite[3.3]{Abu2014-SF} where it was called \emph{normally flat}. The
terminology $e$-flat was first used in \cite{AIKN2018}.

\begin{punto}
Let $F_{S}$ be a right $S$-semimodule. Following Abuhlail \cite{Abu2014-SF},
we say that $F_{S}$ is a \textbf{flat} right $S$-semimodule, if $F$ is the
directed colimit of \emph{finitely presented} projective right $S$%
-semimodules.
\end{punto}

\begin{punto}
Let $M$ be a left $S$-semimodule. A right $S$-semimodules $F$ is called

\textbf{normally }$M$-\textbf{flat}, if for every \emph{subtractive} $S$%
-subsemimodule $L\leq _{S}M,$ we have a \emph{subtractive} submonoid $%
F\otimes _{S}L\leq F\otimes _{S}M;$

$M$-$i$-\textbf{flat}, if for every \emph{subtractive} $S$-subsemimodule $%
L\leq _{S}M,$ we have a submonoid $F\otimes _{S}L\leq F\otimes _{S}M.$

$M$-$m$\textbf{-flat}, if for every $S$-subsemimodule $L\leq _{S}M,$ we have
a submonoid $F\otimes _{S}L\leq F\otimes _{S}M.$

We say that $F_{S}$ is \textbf{normally flat} (resp., $i$-\textbf{flat}, $m$-%
\textbf{flat)}, iff $F$ is normally $M$-flat (resp., $M$-$i$-flat, $M$-$m$%
-flat) for every left $S$-semimodule $M.$
\end{punto}

\begin{defn}
A right $S$-semimodules $F$ is called $M$-$e$\textbf{-flat}, where $M$ is a
left $S$-semimodule, iff for every short exact sequence%
\begin{equation}
0\longrightarrow L\overset{f}{\longrightarrow }M\overset{g}{\longrightarrow }%
N\longrightarrow 0  \label{lmn-33}
\end{equation}%
of left $S$-semimodules the induced sequence%
\begin{equation}
0\longrightarrow F\otimes _{S}L\overset{F\otimes f}{\longrightarrow }%
F\otimes _{S}M\overset{F\otimes g}{\longrightarrow }F\otimes
_{S}N\longrightarrow 0  \label{Ten-33}
\end{equation}%
of commutative monoids is exact. We say that $F_{S}$ is $e$-\textbf{flat},
iff $F$ is $M$-e-flat for every left $S$-semimodule $M$, equivalently the
covariant functor $F\otimes _{S}-:$ $_{S}\mathbf{SM}\longrightarrow $ $_{%
\mathbb{Z}^{+}}\mathbf{SM}$ respects short exact sequences.
\end{defn}

\begin{rem}
The prefix in "$m$-flat" stems from \textbf{mono-flat semimodules}
introduced by Katsov \cite{Kat2004}, and is different from that of $k$%
\textbf{-flat semimodules} in the sense of Al-Thani \cite{Alt2004}, since
the \emph{tensor product} we adopt here is in the sense of Katsov which is
different from that in the sense of Al-Thani (see \cite{Abu2013} for more
details).
\end{rem}

\begin{prop}
\label{unif-flat-e}Let $F$ be a right $S$-semimodule and $M$ a left $S$%
-semimodule.

\begin{enumerate}
\item $F_{S}$ is normally $M$-flat if and only if for every short exact
sequence of the form (\ref{lmn-33}), Sequence (\ref{Ten-33})\ is exact.

\item $F_{S}$ is $M$-$i$-flat if and only if for every short exact sequence
of the form (\ref{lmn-33}), Sequence (\ref{Ten-33})\ is semi-exact, $%
F\otimes f$ is $k$-normal and $F\otimes g$ is normal.

\item $F_{S}$ is $M$-$m$-flat if and only if for every semi-exact sequence
of the form (\ref{lmn-33}) in which $f$ is $k$-normal and $g$ is normal,
Sequence (\ref{Ten-33})\ is semi-exact, $F\otimes f$ is $k$-normal and $%
F\otimes g$ is normal.
\end{enumerate}
\end{prop}

\begin{Beweis}
(1) and (2).

($\Longrightarrow $) Notice that $L=Ker(g)$ is a subtractive $S$%
-subsemimodule of $M.$ Since $F_{S}$ is $M$-$i$-flat (normally $M$-flat), we
know that $F\otimes f$ is a (normal) monomorphism. It follows by Proposition %
\ref{r-exact} (2) (and (3)) that (\ref{Ten-33}) is semi-exact (exact) and $%
F\otimes g$ is a normal epimorphism.

($\Longleftarrow $) Let $L\leq _{S}M$ be a subtractive $S$-subsemimodule.
Then%
\begin{equation}
0\longrightarrow L\overset{\iota }{\longrightarrow }M\overset{\pi _{L}}{%
\longrightarrow }M/L\longrightarrow 0  \label{L/M}
\end{equation}%
is a short exact sequence of left $S$-semimodules, where $\iota $ is the
canonical injection and $\pi _{L}:M\longrightarrow M/L$ is the canonical
projection. By our assumptions, the induced sequence of commutative monoids%
\begin{equation}
0\longrightarrow F\otimes _{S}L\overset{F\otimes \iota }{\longrightarrow }%
F\otimes _{S}M\overset{F\otimes \pi _{L}}{\longrightarrow }F\otimes
_{S}M/L\longrightarrow 0  \label{FotM/L}
\end{equation}%
is semi-exact (exact) and $F\otimes \iota $ is $k$-normal, whence a (normal)
monomorphism.

(3) ($\Longrightarrow $) Since $F_{S}$ is $M$-$m$-flat, we know that $%
F\otimes f$ is a monomorphism, whence $k$-normal. Moreover, it follows by
Proposition \ref{r-exact} (2) that (\ref{Ten-33}) is semi-exact and $%
F\otimes g$ is a normal epimorphism.

($\Longleftarrow $) Let $L\leq _{S}M$ be an $S$-subsemimodule. Then (\ref%
{L/M}) is a semi-exact sequence of left $S$-semimodules in which $\iota $ is
$k$-normal and $\pi _{L}$ is normal. By our assumption, Sequence (\ref%
{FotM/L}) is semi-exact and $F\otimes \iota $ is $k$-normal, whence $%
F\otimes f$ is injective.$\blacksquare $
\end{Beweis}

\begin{prop}
\label{inj->}Let%
\begin{equation}
L\overset{f}{\longrightarrow }M\overset{g}{\longrightarrow }N  \label{LMN-3}
\end{equation}%
be a sequence of left $S$-semimodules, $F$ a right $S$-semimodule and
consider the sequence%
\begin{equation}
F\otimes _{S}L\overset{F\otimes f}{\longrightarrow }F\otimes _{S}M\overset{%
F\otimes g}{\longrightarrow }F\otimes _{S}N  \label{F-LMN3}
\end{equation}%
of commutative monoids.

\begin{enumerate}
\item If (\ref{LMN-3}) is exact with $g$ normal and $F_{S}$ is $e$-flat,
then (\ref{F-LMN3}) is exact and $F\otimes g$ is normal.

\item If (\ref{LMN-3}) is exact with $g$ normal and $F_{S}$ is $i$-flat,
then (\ref{F-LMN3}) is quasi-exact.

\item If (\ref{LMN-3}) is exact and $F_{S}$ is $m$-flat, then (\ref{F-LMN3})
is quasi-exact.
\end{enumerate}
\end{prop}

\begin{Beweis}
By Corollary \ref{M/L}, we have a short exact sequence of left $S$%
-semimodules%
\begin{equation}
0\longrightarrow Ker(g)\overset{\iota }{\longrightarrow }M\overset{\pi }{%
\longrightarrow }M/Ker(g)\longrightarrow 0  \label{Ker-Coker}
\end{equation}%
where $\iota $ and $\pi $ are the canonical $S$-linear maps. Since (\ref%
{LMN-3}) is proper exact, $f(L)=Ker(g)$ and $M/Ker(g)=M/f(L)\simeq Coker(f).$
By the \emph{Universal Property of Kernels}, there exists a unique $S$%
-linear map $\widetilde{f}:L\longrightarrow Ker(g)$ such that $\iota \circ
\widetilde{f}=f$. On the other hand, by the \emph{Universal Property of
Cokernels, }there exists a unique $S$-linear map $\widetilde{g}%
:M/Ker(g)\longrightarrow N$ such that $\widetilde{g}\circ \pi =g$. So, we
have a commutative diagram of left $S$-semimodules%
\begin{equation}
\xymatrix{ & & & L \ar[lld]_{\widetilde{f}} \ar[d]^{f} & & & 0 \ar[ld] \\ 0
\ar[r] & Ker(g) \ar[ld] \ar[rr]^{\iota} & & M \ar[rr]^{\pi} \ar[d]^{g} & &
M/Ker(g) \ar[lld]^{\widetilde{g}} \ar[r] & 0 \\ 0 & & & N & & & }
\label{M-ker}
\end{equation}%
Applying the contravariant functor $F\otimes _{S}-,$ we get the sequence%
\begin{equation}
0\longrightarrow F\otimes _{S}Ker(g)\overset{F\otimes \iota }{%
\longrightarrow }F\otimes _{S}M\overset{F\otimes \pi }{\longrightarrow }%
F\otimes _{S}M/Ker(g)\longrightarrow 0  \label{Ten-Ker}
\end{equation}%
and we obtain the commutative diagram%
\begin{equation}
\xymatrix{& & & F \otimes_S L \ar[dd]^{F \otimes f} \ar[lldd]_{F \otimes
{\widetilde{f}}} & & & 0 \ar[ldd] \\ & & & & & & & \\ 0 \ar[r] & F \otimes_S
Ker(g) \ar[rr]^{F \otimes \iota} \ar[ldd] & & F \otimes_S M \ar[rr]^{F
\otimes \pi} \ar[dd]^{F \otimes g} & & F \otimes_S M/Ker(g) \ar[lldd]^{F
\otimes {\widetilde{g}}} \ar[r] & 0 \\ & & & & & & & \\ 0 & & & F \otimes_S
N & & & &}  \label{Ten-M-ker}
\end{equation}%
of commutative monoids.

Notice that $\widetilde{g}$ is injective since $g$ is $k$-normal. On the
other hand, $\widetilde{f}$ is surjective since $f(L)=Ker(g)$. If $g=%
\widetilde{g}\circ \pi $ is normal, then $\widetilde{g}$ is a normal
monomorphism by Lemma \ref{i-normal} (2-b).

\begin{enumerate}
\item Let $F_{S}$ be $e$-flat and $g$ be normal. Since $F_{S}$ is $e$-flat,
Sequence (\ref{Ten-Ker}) is exact and $F\otimes \widetilde{g}$ is a normal
monomorphism.

\textbf{Step I:}\ We have%
\begin{equation*}
\begin{array}{ccccc}
Ker(F\otimes g) & = & Ker((F\otimes \widetilde{g})\circ (F\otimes \pi )) &
&  \\
& = & Ker(F\otimes \pi ) &  & \text{(}F\otimes \widetilde{g}\text{ is
injective)} \\
& = & im(F\otimes \iota ) &  & \text{((\ref{Ten-Ker}) is (proper-)exact)} \\
& = & im((F\otimes \iota )\circ (F\otimes \widetilde{f})) &  & \text{(}%
F\otimes \widetilde{f}\text{ is surjective)} \\
& = & im(F\otimes f). &  &
\end{array}%
\end{equation*}

\textbf{Step II:}\ Since $F_{S}$ is $e$-flat, it follows by Proposition \ref%
{inj->} (1) that $F\otimes \pi $ is a normal epimorphism. Since $(F\otimes
\widetilde{g})$ is a normal monomorphism, it follows by Lemma \ref{i-normal}
(1-c or 2-c) that $F\otimes g=(F\otimes \widetilde{g})\circ (F\otimes \pi )$
is normal.

\item Let $F_{S}$ be $i$-flat and $g$ be normal. Since $\widetilde{g}$ is a
normal monomorphism and $F_{S}$ be $i$-flat we have: $F\otimes \widetilde{g}$
is a monomorphism. By Proposition \ref{r-exact} (2), Sequence (\ref{Ten-Ker}%
) is quasi-exact, i.e. $\overline{im(F\otimes \iota )}=Ker(F\otimes \pi )$
and $F\otimes \pi $ is a $k$-normal epimorphism. Calculations similar to
those in (1) Step I, show that $\overline{im(F\otimes f)}=Ker(F\otimes g).$
Since $(F\otimes \widetilde{g})$ is a monomorphism, it follows by Lemma \ref%
{i-normal} (1-a) that $F\otimes g=(F\otimes \widetilde{g})\circ (F\otimes
\pi )$ is $k$-normal. Consequently, Sequence (\ref{F-LMN3}) is quasi-exact.

\item Let $F_{S}$ be $m$-flat. Since $g=\widetilde{g}\circ \pi $ is $k$%
-normal, it follows that $\widetilde{g}$ is injective. Since $F_{S}$ is $m$%
-flat, we have: $F\otimes \widetilde{g}$ is a monomorphism. By Proposition %
\ref{r-exact} (3), Sequence (\ref{Ten-Ker}) is quasi-exact, i.e. $\overline{%
im(F\otimes \iota )}=Ker(F\otimes \pi )$ and $F\otimes \pi $ is a $k$-normal
epimorphism. We continue our argument as in (2) to deduce that (\ref{F-LMN3}%
) is quasi-exact.$\blacksquare $
\end{enumerate}
\end{Beweis}

\begin{thm}
\label{e-inj-3}The following are equivalent for a right $S$-semimodule $F:$

\begin{enumerate}
\item $F_{S}$ is normally flat;

\item $F_{S}$ is $e$-flat;

\item For every exact sequence of left $S$-semimodules (\ref{LMN-3}) with $g$
normal, the induced sequence of commutative monoids (\ref{F-LMN3}) is exact
and $F\otimes g$ is normal.
\end{enumerate}
\end{thm}

\begin{Beweis}
$(1)\Longleftrightarrow (2)$ This follows by Proposition \ref{unif-flat-e}
(1).

$(1)\Rightarrow (3)$ This follows by Proposition \ref{inj->} (1).

$(3)\Rightarrow (1)$ This follows directly by applying the assumption to the
exact sequences of left $S$-semimodules of the form $0\longrightarrow M%
\overset{g}{\longrightarrow }N$ with $g$ normal.$\blacksquare $
\end{Beweis}

\begin{thm}
\label{I-inj-3}The following are equivalent for a right $S$-semimodule $F:$

\begin{enumerate}
\item $F_{S}$ is $i$-flat;

\item For every short exact sequence (\ref{lmn-33}) of left $S$-semimodules,
Sequence (\ref{Ten-33}) is semi-exact, $F\otimes f$ is $k$-normal and $%
F\otimes g$ is normal.

\item for every exact\emph{\ }sequence of left $S$-semimodules (\ref{LMN-3})
with $g$ normal, the induced sequence of commutative monoids (\ref{F-LMN3})
is quasi-exact.
\end{enumerate}
\end{thm}

\begin{Beweis}
$(1)\Longleftrightarrow (2)$ This follows by Proposition \ref{unif-flat-e}
(2).

$(1)\Rightarrow (3)$ This follows by Proposition \ref{inj->} (2).

$(3)\Rightarrow (1)$ This follows directly by applying the assumption to the
exact sequences of left $S$-semimodules the form $0\longrightarrow M\overset{%
g}{\longrightarrow }N$ with $g$ normal.$\blacksquare $
\end{Beweis}

\begin{thm}
\label{inj-3}The following are equivalent for a right $S$-semimodule $F:$

\begin{enumerate}
\item $F_{S}$ is $m$-flat;

\item For every semi-exact sequence of the form (\ref{lmn-33}) in which $f$
is $k$-normal and $g$ is normal, Sequence (\ref{Ten-33})\ is semi-exact, $%
F\otimes f$ is $k$-normal and $F\otimes g$ is normal.

\item For every exact sequence (\ref{LMN-3}), the sequence (\ref{F-LMN3}) is
quasi-exact.
\end{enumerate}
\end{thm}

\begin{Beweis}
$(1)\Longleftrightarrow (2)$ This follows by Proposition \ref{unif-flat-e}
(3).

$(1)\Rightarrow (3)$ This follows by Proposition \ref{inj->} (3).

$(3)\Rightarrow (1)$ This follows directly by applying the assumption to the
exact sequences of left $S$-semimodules the form $0\longrightarrow M\overset{%
g}{\longrightarrow }N$.$\blacksquare $
\end{Beweis}

\begin{cor}
\label{S-T}Let $S$ and $T$ be semirings, $F$ a $(T,S)$-bisemimodule and $%
\widetilde{F}$ a right $T$-semimodule. If $F_{S}$ is $e$-flat ($m$-flat) and
$\widetilde{F}_{T}$ is $e$-flat ($m$-flat), then $(\widetilde{F}\otimes
_{T}F)_{S}$ is a normally flat ($m$-flat).
\end{cor}

\begin{Beweis}
Let $F_{S}$ $e$-flat ($m$-flat) and $\widetilde{F}_{T}$ be $e$-flat ($m$%
-flat). By our assumptions and Proposition \ref{unif-flat-e}, the two
functors%
\begin{equation*}
F\otimes _{S}-:\text{ }_{S}\mathbf{SM}\longrightarrow \text{ }_{T}\mathbf{SM}%
\text{ and }\widetilde{F}\otimes _{S}-:\text{ }_{T}\mathbf{SM}%
\longrightarrow \text{ }_{\mathbb{Z}^{+}}\mathbf{SM}
\end{equation*}%
respect short exact sequences (monomorphisms), whence the functor%
\begin{equation*}
\widetilde{F}\otimes _{T}F=(\widetilde{F}\otimes _{S}-)\circ (F\otimes
_{S}-):\text{ }_{S}\mathbf{SM}\longrightarrow \text{ }_{\mathbb{Z}^{+}}%
\mathbf{SM}
\end{equation*}%
respects short exact sequences (monomorphisms). Consequently, $(\widetilde{F}%
\otimes _{T}F)_{S}$ is $e$-flat ($m$-flat).$\blacksquare $
\end{Beweis}

\begin{prop}
\label{flat->e}\emph{(\cite[Theorem 3.6]{Abu2014-SF})} Let $S$ be any
semiring. If $F$ is a flat left (right) $S$-semimodule, then $F$ is $e$-flat.
\end{prop}

\begin{rem}
Let $M$ be a right $S$-semimodule and denote with $\mathcal{F}_{S}^{e}(M)$
(resp. $\mathcal{F}_{S}^{i}(M),$ $\mathcal{F}_{S}^{m}(M)$) the class of $M$-$%
e$-flat (resp. $M$-$i$-flat, $M$-$m$-flat) left $S$-semimodules. Dropping $M$
means that we take the union over all left $S$-semimodules $M.$ Moreover, we
denote with $\mathcal{F}_{S}$ the class of all flat right $S$-semimodules.
It follows directly from the definitions that $M$-$e$-flat and $M$-$m$-flat
right semimodules are $M$-$i$-flat. Moreover, flat semimodules are $m$-flat
by \cite[Proposition 2.1]{Kat2004} $e$-flat by Proposition \ref{flat->e}.
Summarizing, we have the following inclusions:%
\begin{equation}
\mathcal{F}_{S}^{e}(M)\cup \mathcal{F}_{S}^{m}(M)\subseteq \mathcal{F}%
_{S}^{i}(M)\text{ and }\mathcal{F}_{S}\subseteq \mathcal{F}_{S}^{e}\cap
\mathcal{F}_{S}^{m}\subseteq \mathcal{F}_{S}^{i}.  \label{flat-inclusions}
\end{equation}
\end{rem}

\begin{lem}
\label{ret-flat}

\begin{enumerate}
\item Let $M$ be a left $S$-semimodule. Any retract of an $M$-$i$-flat
(resp. $M$-$e$-flat, $M$-$m$-flat) right $S$-semimodule is $M$-$i$-flat
(resp. $M$-$e$-flat, $M$-$m$-flat).

\item Any retract of an $i$-flat (resp. $e$-flat, $m$-flat) right $S$%
-semimodule is $i$-flat (resp. $e$-flat, $m$-flat).
\end{enumerate}
\end{lem}

\begin{Beweis}
We only need to prove \textquotedblleft 1\textquotedblright\ for relative $i$%
-flatness (resp. relative $e$-flatness); the proof for relative $m$-flatness
is similar.

Let $M$ be a left $S$-semimodule, $U\leq _{S}M$ a subtractive subsemimodule,
$F_{S}$ an $M$-$e$-flat right $S$-semimodule and $\widetilde{F}$ a retract
of $F.$ Then there exist $S$-linear maps $\widetilde{F}\overset{\psi }{%
\longrightarrow }F\overset{\theta }{\longrightarrow }\widetilde{F}$ such
that $\theta \circ \psi =\mathrm{id}_{\widetilde{F}}.$ Consider the
commutative diagram%
\begin{equation*}
\xymatrix{0 \ar@{.>}[rr] & & \tilde{F} \otimes_S U
\ar[rr]^{{\tilde{F}}\otimes \iota_U} \ar[d]_{\psi \otimes U} & & {\tilde{F}}
\otimes M \ar[d]^{\psi \otimes M}\\ 0 \ar[rr] & & F \otimes_S U \ar[rr]^{F
\otimes_S \iota_U} \ar[d]_{\theta \otimes U} & & F \otimes_S M
\ar[d]^{\theta \otimes M} \\ & & {\tilde{F}} \otimes_S U
\ar[rr]^{{\tilde{F}}\otimes \iota_U} & & {\tilde{F}} \otimes_S M}
\end{equation*}%
Indeed, $(\theta \otimes _{S}\mathrm{id}_{U})\circ (\psi \otimes _{S}\mathrm{%
id}_{U})=\mathrm{id}_{\widetilde{F}\otimes _{S}U}$ and $(\theta \otimes _{S}%
\mathrm{id}_{M})\circ (\psi \otimes _{S}\mathrm{id}_{M})=\mathrm{id}_{%
\widetilde{F}\otimes _{S}M},$ \emph{i.e.} $\widetilde{F}\otimes _{S}U$ is a
retract of $F\otimes _{S}U$ and $\widetilde{F}\otimes _{S}M$ is a retract of
$F\otimes _{S}M.$ In particular, $\psi \otimes U$ and $\psi \otimes M$ are
monomorphisms.

If $F_{S}$ is $M$-$i$-flat (resp. $M$-$e$-flat), then $\mathrm{id}%
_{F}\otimes _{S}\iota _{U}:F\otimes _{S}U\longrightarrow F\otimes _{S}M$ is
(normal) monomorphism. It follows that $\mathrm{id}_{\widetilde{F}}\otimes
_{S}\iota _{U}$ is a (normal) monomorphism by Lemma \ref{i-normal}
\textquotedblleft 1\textquotedblright , \emph{i.e.} $\widetilde{F}\otimes
_{S}U\leq _{S}\widetilde{F}\otimes _{S}M$ is a (subtractive) $S$-semimodule.
Consequently, $\widetilde{F}$ is $M$-$i$-flat ( $M$-$e$-flat).$\blacksquare $
\end{Beweis}

\begin{prop}
\label{sum-flat}Let $\{F_{\lambda }\}_{\Lambda }$ be a non-empty collection
of right $S$-semimodules.

\begin{enumerate}
\item Let $M$ be a left $S$-semimodule. Then $\bigoplus\limits_{\lambda \in
\Lambda }F_{\lambda }$ is $M$-$i$-flat (resp. $M$-$e$-flat, $M$-$m$-flat) if
and only if $F_{\lambda }$ is $M$-$i$-flat (resp. $M$-$e$-flat, $M$-$m$%
-flat) for every $\lambda \in \Lambda .$

\item $\bigoplus\limits_{\lambda \in \Lambda }F_{\lambda }$ is $i$-flat
(resp. $e$-flat, $m$-flat) if and only if $F_{\lambda }$ is $i$-flat (resp. $%
e$-flat, $m$-flat) for every $\lambda \in \Lambda .$
\end{enumerate}
\end{prop}

\begin{Beweis}
We only need to prove \textquotedblleft 1\textquotedblright\ for relative $i$%
-flatness (resp. relative $e$-flatness); the proof for relative $m$-flatness
is similar (cf. \cite[Proposition 2.3]{Alt2004} for $k$-flat semimodules).

($\Longrightarrow $)\ For every $\lambda \in \Lambda ,$ $F_{\lambda }$ is a
retract of\emph{\ }$\bigoplus\limits_{\lambda \in \Lambda }F_{\lambda },$
whence $M$-$i$-flat (resp. $M$-$e$-flat, $M$-$m$-flat) by Lemma \ref%
{ret-flat}.

($\Longleftarrow $)\ Let $F:=$ $\bigoplus\limits_{\lambda \in \Lambda
}F_{\lambda }$ and consider the projections $\pi _{\lambda
}:F\longrightarrow F_{\lambda },$ $(f_{\lambda })_{\Lambda }\mapsto
f_{\lambda }$ for $\lambda \in \Lambda .$ Let $U\leq _{S}M$ be a subtractive
$S$-subsemimodule. Assume that $F_{\lambda }$ is $M$-$i$-flat ($M$-$e$-flat)
for every $\lambda \in \Lambda .$ Then $F_{\lambda }\otimes _{S}U\leq
_{S}F_{\lambda }\otimes _{S}M$ is a (subtractive) subsemimodule for every $%
\lambda \in \Lambda ,$ whence $\bigoplus\limits_{\lambda \in \Lambda
}(F_{\lambda }\otimes _{S}U)\leq _{S}\bigoplus\limits_{\lambda \in \Lambda
}(F_{\lambda }\otimes _{S}M)$ is a (subtractive) subsemimodule by Lemma \ref%
{u-sum} (1).

Since%
\begin{equation*}
\bigoplus\limits_{\lambda \in \Lambda }(F_{\lambda }\otimes _{S}U)\simeq
\bigoplus\limits_{\lambda \in \Lambda }F_{\lambda }\otimes _{S}U\text{ and }%
\bigoplus\limits_{\lambda \in \Lambda }(F_{\lambda }\otimes _{S}M)\simeq
\bigoplus\limits_{\lambda \in \Lambda }F_{\lambda }\otimes _{S}M,
\end{equation*}%
we conclude that $\bigoplus\limits_{\lambda \in \Lambda }F_{\lambda }\otimes
_{S}U\leq _{S}\bigoplus\limits_{\lambda \in \Lambda }F_{\lambda }\otimes
_{S}M$ is a (subtractive) subsemimodule. It follows that $%
\bigoplus\limits_{\lambda \in \Lambda }F_{\lambda }$ is $M$-$i$-flat ($M$-$e$%
-flat).$\blacksquare $
\end{Beweis}

\begin{prop}
\label{closed-sub-factor}Let $F$ be a right $S$-semimodule and%
\begin{equation*}
0\longrightarrow L\overset{f}{\longrightarrow }M\overset{g}{\longrightarrow }%
N\longrightarrow 0
\end{equation*}%
be an exact sequence of left $S$-semimodules.

\begin{enumerate}
\item If $F$ is $M$-$i$-flat (resp. $M$-$e$-flat, $M$-$m$-flat), then $F$ is
$L$-$i$-flat (resp. $L$-$e$-flat, $L$-$m$-flat).

\item If $F$ is $M$-$m$-flat (resp. $M$-$i$-flat), then $F$ is $N$-$m$-flat
(resp. $N$-$i$-flat).
\end{enumerate}
\end{prop}

\begin{Beweis}
\begin{enumerate}
\item Let $F_{S}$ be $M$-$i$-flat (resp. $M$-$e$-flat) and $U\leq _{S}L$ be
a subtractive $S$-subsemimodule. Then $U\leq _{S}M$ is a subtractive $S$%
-semimodule (notice that $L\leq _{S}M$ is subtractive). Since $F$ is $M$-$i$%
-flat ($M$-$e$-flat), $F\otimes _{S}U\leq F\otimes _{S}M$ and $F\otimes
_{S}L\leq F\otimes _{S}M$ are (subtractive) submonoids, whence $F\otimes
_{S}U\leq F\otimes _{S}L$ is a (subtractive) submonoid (by Lemma \ref%
{i-normal} (1-b)). Consequently, $F$ is $L$-$i$-flat ($L$-$e$-flat).
Similarly, one can prove that if $F_{S}$ is $M$-$m$-flat, then $F_{S}$ is $L$%
-$m$-flat.

\item Consider the \emph{pullback} $(P,\iota ^{\prime },g^{\prime })$ of $%
\iota :U\hookrightarrow N$ and $g:M\longrightarrow N$ given by Lemma \ref%
{pull}. Applying $F\otimes _{S}-$ to Diagram (\ref{pbc}) yields the
following commutative diagram%
\begin{equation*}
\xymatrix{& & 0 \ar[d] & 0 \ar@{.>}[d] & \\ & F\otimes_S L \ar[r]^{F\otimes
f'} \ar@{=}[d]_{F\otimes id} & F\otimes_S P \ar[d]^{F\otimes g'}
\ar[r]^{F\otimes \iota '} & F\otimes_SU \ar[d]^{F\otimes \iota} \ar[r] & 0
\\ & F\otimes_S L \ar[r]_{F\otimes f} & F\otimes_S M \ar[r]_{F\otimes g} &
F\otimes_S N \ar[r] & 0}
\end{equation*}%
of commutative monoids. Assume that $F_{S}$ is $M$-$m$-flat (resp. $M$-$i$%
-flat), so that $F\otimes g^{\prime }$ is injective. Since the rows are
semi-exact, $F\otimes \iota ^{\prime }$ is surjective and $F\otimes g$ is a
normal epimorphism (by Proposition \ref{r-exact}), it follows by Lemma \ref%
{semi-ex} (1) that $F\otimes \iota $ is injective.$\blacksquare $
\end{enumerate}
\end{Beweis}

\begin{lem}
\label{ds}Let $F$ be a right $S$-semimodule.

\begin{enumerate}
\item Let $M$ be a left $S$-semimodule. Then $F$ is $M$-$m$-flat if and only
if for every finitely generated $S$-subsemimodule and exact sequence $%
0\longrightarrow K\overset{\iota _{K}}{\longrightarrow }M$ the sequence of
commutative monoid $0\longrightarrow F\otimes _{S}K\overset{F\otimes \iota
_{K}}{\longrightarrow }F\otimes _{S}M$ is exact.

\item Let $L$ and $N$ be cancellative left $S$-semimodules. If $F$ is \emph{%
cancellative}, $L$-$m$-flat and $N$-$m$-flat, then $F$ is $L\oplus N$-$m$%
-flat.

\item Let $\{M_{\lambda }\}_{\Lambda }$ be a collection of \emph{cancellative%
} left $S$-semimodules. If $F$ is cancellative and $M_{\lambda }$-$m$-flat
for every $\lambda \in \Lambda ,$ then $F$ is $\bigoplus\limits_{\lambda \in
\Lambda }M_{\lambda }$-$m$-flat.
\end{enumerate}
\end{lem}

\begin{Beweis}
\begin{enumerate}
\item $(\Longrightarrow )$ Obvious.

$(\Longleftarrow )$ Let $U\leq _{S}M.$ Suppose that
\begin{equation*}
(F\otimes \iota _{U})(\sum\limits_{i=1}^{m}f_{i}\otimes _{S}u_{i})=(F\otimes
\iota _{U})(\sum\limits_{j=1}^{n}\widetilde{f}_{j}\otimes _{S}\widetilde{u}%
_{j}),
\end{equation*}%
and consider $K\leq _{S}M$ generated by $\{u_{1},\cdots ,u_{m},\widetilde{u}%
_{1},\cdots ,\widetilde{u}_{n}\}.$ Notice that%
\begin{eqnarray*}
(F\otimes \iota _{K})(\sum\limits_{i=1}^{m}f_{i}\otimes _{S}u_{i})
&=&(F\otimes \iota _{U})(\sum\limits_{i=1}^{m}f_{i}\otimes _{S}u_{i}) \\
&=&(F\otimes \iota _{U})(\sum\limits_{j=1}^{n}\widetilde{f}_{j}\otimes _{S}%
\widetilde{u}_{j}) \\
&=&(F\otimes \iota _{K})(\sum\limits_{j=1}^{n}\widetilde{f}_{j}\otimes _{S}%
\widetilde{u}_{j}),
\end{eqnarray*}%
whence $\sum\limits_{i=1}^{m}f_{i}\otimes _{S}u_{i}=\sum\limits_{j=1}^{n}%
\widetilde{f}_{j}\otimes _{S}\widetilde{u}_{j}$ since $F\otimes \iota _{K}$
is injective. It follow that $F\otimes \iota _{U}$ is injective.
Consequently, $F$ is $M$-$m$-flat.

\item Let $U\leq _{S}L\oplus N$ and consider the short exact sequence%
\begin{equation*}
0\longrightarrow L\overset{\iota }{\longrightarrow }L\oplus N\overset{\pi }{%
\longrightarrow }N\longrightarrow 0
\end{equation*}%
of cancellative left $S$-semimodules. Consider the pullback $(P,\lambda
^{\prime },\iota ^{\prime })$ of $\lambda :U\hookrightarrow L\oplus N$ and $%
\iota :L\hookrightarrow L\oplus N$ given by%
\begin{eqnarray}
P &=&\{(u,l)\in U\times L\text{ }|\text{ }\lambda (u)=\iota (l)\};
\label{PUL} \\
\lambda ^{\prime } &:&P\rightarrow U,\text{ }(u,l)\mapsto u;  \notag \\
\iota ^{\prime } &:&P\rightarrow L,\text{ }(u,l)\mapsto l  \notag
\end{eqnarray}%
and the commutative diagram%
\begin{equation}
\xymatrix{& 0 \ar[d] & 0 \ar[d] & 0 \ar[d] & \\ 0 \ar[r] & P \ar[r]^{\lambda
'} \ar_{\iota '}[d] & U \ar[d]^{\lambda} \ar[r]^{\pi '} & U/P \ar[d]^{h}
\ar[r] & 0 \\ 0 \ar[r] & L \ar[r]_{\iota} & L \oplus N \ar[r]_{\pi} & N
\ar[r] & 0}  \label{diag-pul}
\end{equation}

of cancellative $S$-semimodules. Applying $F\otimes _{S}-$ to Diagram (\ref%
{diag-pul}) yields the following commutative diagram%
\begin{equation}
\xymatrix{ & 0 \ar[d] & 0 \ar@{.>}[d] & 0 \ar[d] & \\ & F \otimes _S P
\ar[r]^{F \otimes \lambda '} \ar_{F \otimes \iota '}[d] & F \otimes_S U
\ar[d]^{F \otimes \lambda} \ar[r]^{F \otimes \pi '} & F \otimes_S U/P
\ar[d]^{F \otimes h} \ar[r] & 0 \\0 \ar[r] & F \otimes L \ar[r]_{F \otimes
\iota} & F \otimes_S (L \oplus N) \ar[r]_{F \otimes \pi} & F \otimes_S N
\ar[r] & 0}  \label{diag-F}
\end{equation}%
of cancellative commutative monoids in which the second row is exact. By
Proposition \ref{r-exact}, the first row is semi-exact and $F\otimes \pi
^{\prime }$ is a normal epimorphism. Since $F$ is $L$-$m$-flat and $N$-$m$%
-flat, both $F\otimes \iota ^{\prime }$ and $F\otimes h$ are injective. It
follows by Lemma \ref{semi-ex} (2) that $F\otimes \lambda $ is injective.
Consequently, $F$ is $L\oplus N$-$m$-flat.

\item Let $U\leq _{S}\bigoplus\limits_{\lambda \in \Lambda }M_{\lambda }.$
In light of (1), we can assume that $U$ is finitely generated, whence
contained in a finite number of direct sums. So, we are done by (2).$%
\blacksquare $
\end{enumerate}
\end{Beweis}

\section{Von Neumann Regular Semirings}

\qquad In this section, we study the so called von Neumann regular semirings
that are not necessarily rings.

\begin{defn}
A semiring $S$ is

a \textbf{von Neumann regular} semiring, iff for every $a\in S$ there exists
some $s\in S$ such that $a=asa;$

an \textbf{additively-regular semiring}, iff for every $a\in S,$ there
exists $b\in S$ such that $a+b+a=a.$
\end{defn}

It is well known that for a ring $R$ and a positive integer $n,$ the ring $R$
is von Neumann Regular if and only if the matrix ring $M_{n}(R)$ is von
Neumann regular. The situation for semirings is different (e.g., \cite[%
Theorem 2.24]{Kap1969}). While the von Neumann regularity of a matrix
semiring $M_{n}(S)$ over a semiring $S$ implies indeed that $S$ is von
Neumann regular, the converse is not true in general. Moreover, for $n\geq
3, $ $M_{n}(S)$ is a von Neumann regular if and only if $S$ is a von Neumann
regular ring \cite[Proposition 2]{Ili2001} (cf. \cite[page 70, Proposition
5.25]{Gol2003}, \cite[Theorem 2.3.]{CSL2009}, \cite[Proof of Proposition 4.9]%
{AIKN2015}).

\begin{ex}
(\cite[Remark 2.2.]{CSL2009}) Let $S=(\{0,1,2,3\},\max ,0,\min ,3).$ One can
easily check that $S$ is von Neumann regular, however $A=\left[
\begin{array}{cc}
0 & 1 \\
2 & 3%
\end{array}%
\right] \in M_{2}(S)$ is \emph{not }a regular element (i.e. there is no $%
B\in M_{2}(S)$ with $ABA=A.$
\end{ex}

Assuming all semimodules of a given \emph{commutative} semiring $S$ to be
(mono-)flat forces the semiring to be a von Neumann regular \emph{ring}
(cf., \cite[Theorem 2.11]{Kat2004}. This suggests other notions of flatness,
e.g. $e$-flatness and $i$-flatness.

\begin{definition}
\cite[page 71]{Gol1999} Let $S$ be a semiring. We say that $S$ is a \textbf{%
left subtractive semiring }(\textbf{right subtractive semiring})%
\index{Semiring!subtractive} if every left (right) ideal of $S$ is
subtractive. We say that $S$ is a subtractive semiring if $S$ is both left
and right subtractive. 
\end{definition}

\begin{rem}
Whether a left subtractive semiring is necessarily right subtractive was an
open problem till a counterexample was given in \cite[Fact 2.1]{KNT2011}.
\end{rem}

\subsection*{Homological Results}

Von Neumann regular rings are characterized by the fact that all left
(right) modules over them are flat, whence called \emph{absolutely flat rings%
} by the Bourbakian school. We generalize this result partially by showing
that a sufficient condition for a (left and right) subtractive semiring $S$
to be von Neumann regular is the assumption that all left (right) $S$%
-semimodules are $S$-$e$-flat.

\bigskip

The proof of the following lemma is obtained by \emph{diagram chasing}, with
appropriate modifications, which is a well-known tool in the classical
proofs which can be found in standard book of Homological Algebra (cf., \cite%
[Proposition 2.70, Corollary 3.59, Proposition 3.60]{Rot2009}).

\begin{lem}
\label{lem359}Let $A$ be a right $S$-semimodule and consider for every left
ideal $I\leq _{S}S$ the canonical surjective map of commutative monoids%
\begin{equation}
\theta _{I}:A\otimes _{S}I\longrightarrow AI,%
\text{ }a\otimes _{S}i\mapsto ai.  \label{th_A}
\end{equation}

\begin{enumerate}
\item $A_{S}$ is $S$-$m$-flat if and only if $A\otimes _{S}I\overset{\theta
_{I}}{\simeq }AI$ for every (finitely generated)\ left ideal $I$ of $S.$

\item $A_{S}$ is $S$-$i$-flat if and only if $A\otimes _{S}I\overset{\theta
_{I}}{\simeq }AI$ for every subtractive left ideal $I$ of $S.$

\item $A_{S}$ is $S$-$e$-flat if and only if $A\otimes _{S}I\overset{\theta
_{I}}{\simeq }AI$ and $AI\leq A$ is subtractive for every subtractive left
ideal $I$ of $S.$
\end{enumerate}
\end{lem}

\begin{Beweis}
Let $I$ be a left ideal of $S.$ Consider the embeddings $\iota
_{I}:I\hookrightarrow S,$ $\zeta _{I}:AI\hookrightarrow A$ and recall the
canonical isomorphism $A\otimes _{S}S\overset{\varphi _{A}}{\simeq }A$
(Lemma \ref{lem258}). Then we have the commutative diagram of commutative
monoids%
\begin{equation*}
\xymatrix{& & 0 \ar[d] \\ & A \otimes_S I \ar[r]^{A \otimes \iota_I}
\ar[d]_{\theta_I} & A \otimes_S S \ar[d]^{\varphi _{A}} \\ 0 \ar[r] & AI
\ar[r]_{\zeta _{I}} & A}
\end{equation*}%
The result follows now directly from the definitions noticing that for every
left ideal of $A$ we have $A\otimes \iota _{I}$ is injective if and only if $%
\theta _{I}$ is injective. In light of Lemma \ref{ds}, it is sufficient in
(1) to consider only the finitely generated left ideals of $S$.$\blacksquare
$
\end{Beweis}

\begin{rem}
Part ($\Longrightarrow $) of (3) in Lemma \ref{lem359} was proved in \cite[%
Theorem 2, Corollary 1]{NG2019} assuming the semiring $S$ is commutative.
Notice that (2) provides a complete characterization of right $S$%
-semimodules $A_{S}$ for which $A\otimes _{S}I\overset{\theta _{I}}{\simeq }%
AI$ is an isomorphism for every subtractive left ideal $I$ of $S.$
\end{rem}

\begin{defn}
We say that a left $S$-semimodule $M$ is \textbf{normally }$S$\textbf{%
-generated}, if there exists a \emph{normal} epimorphism $S^{(\Lambda )}%
\overset{\pi }{\longrightarrow }M\longrightarrow 0.$ We say that $_{S}S$ is
a \textbf{normal generator} iff every left $S$-semimodule is normally $S$%
-generated.
\end{defn}

\begin{prop}
Let $S$ be a cancellative semiring and $F$ a cancellative right $S$%
-semimodule. The following are equivalent:

\begin{enumerate}
\item $F_{S}$ is $S$-$m$-flat;

\item The canonical map $\theta _{I}:F\otimes _{S}I\longrightarrow FI$ of
commutative monoids is injective, whence an isomorphism, for every (finitely
generated) left ideal $I$ of $S;$

\item $F_{S}$ is $N$-$m$-flat for every normally $S$-generated left $S$%
-semimodule $N.$
\end{enumerate}
\end{prop}

\begin{Beweis}
The equivalence $(1)\Longleftrightarrow (2)$ follows by Lemma \ref{lem359}
(without assuming that $S$ is cancellative). The implication $(3)\Rightarrow
(1)$ is trivial. Assume $(1).$ Let $N$ be normally $S$-generated so that
there exists a normal epimorphism $\pi :S^{(\Lambda )}\longrightarrow N$ for
some index set $\Lambda .$ Consider the short exact sequence%
\begin{equation*}
0\longrightarrow Ker(\pi )\longrightarrow S^{(\Lambda )}\overset{\pi }{%
\longrightarrow }N\longrightarrow 0.
\end{equation*}%
Since $F_{S}$ is $S$-$m$-flat by (1), it follows by Lemma \ref{ds} that $%
F_{S}$ is $S^{(\Lambda )}$-$m$-flat. Then $F_{S}$ is $N$-flat by Proposition %
\ref{closed-sub-factor}.$\blacksquare $
\end{Beweis}

The assumptions of the following result hold in particular when $S$ is a
ring, whence it recovers the classical result (e.g., \cite[12.6]{Wis1991}).

\begin{cor}
\label{S-mon}Let $S$ be a cancellative semiring such that $_{S}S$ is a
normal generator. A cancellative right $S$-semimodule $F$ is $m$-flat if and
only if $F_{S}$ is $S$-$m$-flat.
\end{cor}

\begin{lem}
\label{lem360}Let $F$ be an $S$-$m$-flat ($S$-$i$-flat)\ right $S$%
-semimodule and $K\overset{\iota }{\hookrightarrow }F$ a subtractive $S$%
-subsemimodule.

\begin{enumerate}
\item If $F/K$ is $S$-$m$-flat ($S$-$i$-flat) and $KI\leq _{S}K$ is
subtractive, then $K\cap FI=KI$ for every (subtractive) left ideal $I$ of $%
S. $

\item If $K\cap FI=KI$ for every finitely generated left ideal of $S,$ then $%
F/K$ is $S$-$m$-flat.

\item If $K\cap FI=KI$ (and $FI\leq F$ is subtractive)\ for every
subtractive left ideal $I$ of $S,$ then $F/K$ is $S$-$i$-flat ($S$-$e$-flat).
\end{enumerate}
\end{lem}


\begin{Beweis}
Consider the right $S$-semimodule $A:=F/K$ and recall, by Lemma \ref{exact}
(7), that we have a short exact sequence of right $S$-semimodules%
\begin{equation}
0\rightarrow K\overset{\iota }{\longrightarrow }F\overset{\varphi }{%
\longrightarrow }A\rightarrow 0.  \label{KFA}
\end{equation}

Let $I\leq _{S}S$ be an arbitrary (subtractive) left ideal. Applying $%
-\otimes _{S}I$ to the exact sequence (\ref{KFA}), it follows by Lemma \ref%
{r-exact} (3) that the following sequence%
\begin{equation*}
K\otimes _{S}I\overset{\iota \otimes I}{\longrightarrow }F\otimes _{S}I%
\overset{\varphi \otimes I}{\longrightarrow }A\otimes _{S}I\rightarrow 0
\end{equation*}
of commutative monoids is semi-exact and $\varphi \otimes I$ is a normal
epimorphism. Consider the following commutative diagram%
\begin{equation}
\xymatrix{K \otimes_S I \ar[r]^{\iota\otimes I} \ar[d]_{\theta_K} & F
\otimes_S I \ar[r]^{\varphi\otimes I} \ar[d]_{\theta_F} & A \otimes_S I
\ar@{-->}[d]^{\gamma} \ar[r] & 0\\ KI \ar[r]^{\iota'} & FI \ar[r]^{\pi} &
FI/KI \ar[r] & 0}  \label{diag-KF}
\end{equation}%
of commutative monoids with semi-exact rows.

Notice that $\theta _{F}$ is injective, whence an isomorphism, since $F_{S}$
is $S$-$m$-flat ($S$-$i$-flat). Since $\varphi \otimes I$ and $\pi $ are
normal epimorphisms, $\theta _{K}$ is surjective and $\theta _{F}$ is
injective, there exists by Lemma \ref{semi-ex} a unique isomorphism
\begin{equation*}
\gamma :A\otimes _{S}I\longrightarrow FI/KI
\end{equation*}%
of commutative monoids that makes Diagram (\ref{diag-KF}) commute.

Since $\varphi :F\longrightarrow A$ is surjective, $\varphi (FI)=AI.$
Consider the restriction $\varphi |_{FI}:FI\rightarrow AI$ and notice that $%
Ker(\varphi |_{FI})=FI\cap K.$ Consider
\begin{equation*}
\beta :AI\rightarrow FI/(FI\cap K),\text{ }ai\mapsto \lbrack fi]\text{ where
}\varphi (f)=a.
\end{equation*}

\vspace{0.5cm}

\textbf{Claim I:\ }$\beta $ is well-defined.

Suppose that $\varphi (f)=a=\varphi (f^{\prime })$ for some $f,f^{\prime
}\in F.$ Since $\varphi $ is $k$-normal, there exist $k,k^{\prime }\in K$
such that $f+k=f^{\prime }+k^{\prime },$ whence $fi+ki=f^{\prime
}i+k^{\prime }i$ for every $i\in I.$ Since $ki,k^{\prime }i\in FI\cap K$, we
get $[fi]=[f^{\prime }i]$. So, $\beta $ is well defined as it is well
defined on a generating set of $AI$.

\vspace{0.5cm}

\textbf{Claim II: }$\beta $ is injective.

Suppose that $\beta (\sum\limits_{j}a_{j}i_{j})=\beta
(\sum\limits_{j}a_{j}^{\prime }i_{j}^{\prime })$ for some $%
\sum\limits_{j}a_{j}i_{j},$ $\sum\limits_{j}a_{j}^{\prime }i_{j}^{\prime }$ $%
\in AI.$ Then $[\sum\limits_{j}f_{j}i_{j}]=[\sum\limits_{j}f_{j}^{\prime
}i_{j}^{\prime }]$ for some $f_{j},f_{j}^{\prime }\in F$ satisfying $\varphi
(f_{j})=a_{j}$ and $\varphi (f_{j}^{\prime })=a_{j}^{\prime }$. It follows
that $\sum\limits_{j}f_{j}i_{j}+z=\sum\limits_{j}f_{j}^{\prime
}i_{j}^{\prime }+z^{\prime }$ for some $z,z^{\prime }\in FI\cap K$ and so%
\begin{equation*}
\begin{array}{ccccc}
\sum\limits_{j}a_{j}i_{j} & = & \sum\limits_{j}\varphi (f_{j})i_{j} & = &
\varphi (\sum\limits_{j}f_{j}i_{j}+z) \\
& = & \varphi (\sum\limits_{j}f_{j}^{\prime }i_{j}^{\prime }+z^{\prime }) & =
& \sum\limits_{j}\varphi (f_{j}^{\prime })i_{j}^{\prime } \\
& = & \sum\limits_{j}a_{j}^{\prime }i_{j}^{\prime }. &  &
\end{array}%
\end{equation*}

Consider the commutative diagram
\begin{equation*}
\xymatrix{ & & A \otimes_S I \ar[rr]^{\gamma} \ar[dd]_{\theta_A} & & FI/KI
\ar[dd]^{\sigma} \\ & & & & \\ 0 \ar[rr] & & AI \ar[rr]_{\beta} \ar[d] & &
FI/(FI \cap K) \ar[d] \\ & & 0 & & 0 }
\end{equation*}%
where%
\begin{equation*}
\sigma :FI/KI\longrightarrow FI/(FI\cap K),\text{ }[fi]_{KI}\mapsto \lbrack
fi]_{FI\cap K}.
\end{equation*}%
Since $\gamma $ is an isomorphism, we conclude that $\sigma $ is injective
(whence an isomorphism) if and only if $\theta _{A}$ is injective (an
isomorphism).

\begin{enumerate}
\item Let $A$ be $S$-$m$-flat ($S$-$i$-flat). In this case, $\theta _{A}$ is
an isomorphism for every (subtractive) left ideal $I\leq _{S}S$ by Lemma \ref%
{lem359} and it follows that $\sigma $ is injective. In particular, $(FI\cap
K)/KI=Ker(\sigma )=0.$ Since $KI\leq _{S}K$ is subtractive (by assumption),
we conclude that $KI=FI\cap K.$

\item If $FI\cap K=KI$ for any finitely generated left ideal $I$ of $S,$
then $\sigma $ is injective, whence $\theta _{A}$ is injective. The result
follows now by Lemma \ref{lem359}.

\item The proof is similar to that of (2).$\blacksquare $
\end{enumerate}
\end{Beweis}

The proof of the following lemma is similar to that of von Neumann regular
rings (e.g., the equivalence of (a) and (b) in \cite[3.10]{Wis1991}).

\begin{lem}
\label{regular-idemp}A semiring $S$ is von Neumann regular semiring if and
only if every principal left (right)\ ideal of $S$ is generated by an
idempotent.
\end{lem}

Von Neumann regular rings are characterized by the fact that each principal
(finitely generated) left (right) ideal is a direct summand \cite[3.10]%
{Wis1991}. The following \emph{counterexample}, communicated to the Authors
by T. Nam, shows that this property is not necessarily true for von Neumann
regular proper semirings.

\begin{ex}
Consider a chain $0<a<1.$ Then $S:=(\{0,a,1\},\max ,0,\min ,1)$ is an
(additively-regular) von Neumann regular commutative semiring. The principal
ideal $I:=Sa$ is \emph{not} a direct summand of $S.$
\end{ex}

We call semiring $S$ \emph{left}\ (\emph{right})\ \emph{B\'{e}zout}, iff
every finitely generated left (right) ideal is principal. We call $S$ a B%
\'{e}zout semiring, iff $S$ is left B\'{e}zout and right B\'{e}zout.

\begin{prop}
\label{Bez-Neumann}If $S$ is a left (right) B\'{e}zout von Neumann regular
semiring, then every normally $S$-generated right (left)$\ S$-semimodule is $%
S$-$m$-flat.
\end{prop}

\begin{Beweis}
Assume that $S$ is a left B\'{e}zout von Neumann regular semiring. Let $A$
be a normally $S$-generated right $S$-semimodule. Then there exists an exact
sequence of left $S$-semimodules%
\begin{equation*}
0\longrightarrow K\overset{\iota }{\longrightarrow }F\overset{\pi }{%
\longrightarrow }A\longrightarrow 0
\end{equation*}%
where $F\simeq S^{(\Lambda )}$ for some index set $\Lambda ,$ and $%
K:=Ker(\pi ).$ Since $F_{S}$ is free, it is flat and in particular $m$-flat.
Let $I$ be a finitely generated left ideal of $S.$ Since $S$ is left B\'{e}%
zout, it $I$ is principal and so $I=Se$ for some idempotent $e$ of $S$ (by
Lemma \ref{regular-idemp}). Since $S$ is von Neumann regular, there exists
some $e^{\prime }\in S$ such that $e=ee^{\prime }e.$ Let $k=fe\in FI\cap K$
for some $k\in K$ and $f\in F.$ Then%
\begin{equation*}
k=fe=f(ee^{\prime }e)=(fe)(e^{\prime }e)=(ke^{\prime })e\in KI.
\end{equation*}%
The result follows now by Lemma \ref{lem360} (2).$\blacksquare $
\end{Beweis}

\begin{ex}
\label{ABC}Let $S$ be an additively-regular von Neumann regular semiring.
For each $a\in S,$ denote by $a^{\prime }\in S$ the \emph{unique}\textit{\ }%
element in $S$ satisfying $a+a^{\prime }+a=a$ and $a^{\prime }+a+a^{\prime
}=a^{\prime }$ (cf., \cite[Proposition 13.1]{Gol1999}). If all elements $%
a,b\in S$ satisfy
\begin{equation*}
\text{(A) }a(a+a^{\prime })=a+a^{\prime };\text{ (B) }a(b+b^{\prime
})=(b+b^{\prime })a;\text{ (C) }a+a(b+b^{\prime })=a,
\end{equation*}%
then $S$ is left (and right)\ B\'{e}zout by \cite{SM2006}.
\end{ex}

\qquad As a direct consequence of Proposition \ref{Bez-Neumann} and Example %
\ref{ABC}, we obtain:

\begin{cor}
Let $R$ be a von Neumann regular ring, $D$ a distributive lattice and
consider the additively-regular von Neumann regular semiring $S:=R\times D.$
Then every normally $S$-generated right (left)$\ S$-semimodule is $S$-$m$%
-flat.
\end{cor}

\begin{thm}
\label{sflatvon}If $S$ is subtractive and every right (left) $S$-semimodule
is $S$-$e$-flat, then $S$ is a von Neumann regular semiring.
\end{thm}

\begin{Beweis}
Let $a\in S.$ By our assumption, $S$ is \emph{right subtractive}, whence $%
K:=aS$ is a subtractive right ideal of $S$ and%
\begin{equation*}
0\longrightarrow aS\longrightarrow S\longrightarrow S/aS\rightarrow 0
\end{equation*}%
is an exact sequence of right $S$-semimodules by Lemma \ref{exact} (7).
Indeed, $F:=S_{S}$ is ($S$)-$e$-flat. By our assumptions, the right $S$%
-semimodules $aS$ and $S/aS$ are both $S$-$e$-flat and so it follows, by
Lemma \ref{lem360} (1), that for every \emph{subtractive} left ideal $I$ of $%
S:$%
\begin{equation*}
aS\cap I=aS\cap SI=K\cap FI=KI=(aS)I.
\end{equation*}%
By our assumption, $S$ is \emph{left subtractive} and so the left ideal $%
I:=Sa\leq _{S}S$ is subtractive, whence%
\begin{equation*}
aS\cap Sa=(aS)(Sa)=aSa.
\end{equation*}%
It follows that $a\in aSa,$ \emph{i.e.} exists some $s\in S$ such that $%
a=asa.\blacksquare $
\end{Beweis}

\begin{cor}
\label{flatvon}If $S$ is subtractive commutative semiring such that every $S$%
-semimodule is $S$-$e$-flat, then $S$ is a von Neumann regular semiring.
\end{cor}

\bigskip

In light of Theorem \ref{sflatvon} and the fact that a commutative semiring
over which all semimodules are flat is a von Neumann regular \emph{ring}
(cf., \cite[Theorem 2.11.]{Kat2004}) we raise the following question:

\bigskip

\textbf{Question:\ }Which class of subtractive semirings is characterized by
the $e$\emph{-flatness} of all right (left) semimodules?


\end{document}